\documentclass[reqno,12pt]{amsart}

\usepackage[all]{xy}
\usepackage{amssymb}
\usepackage{amsgen}
\usepackage{amsmath}
\usepackage{amsthm}
\usepackage{mathrsfs}
\usepackage{cite}
\usepackage{amsfonts}

\renewcommand{\to}{\longrightarrow}

\newcommand{\J}{{\mathrel{\mathscr J}}} 
\newcommand{\D}{{\mathrel{\mathscr D}}} 
\newcommand{\R}{{\mathrel{\mathscr R}}} 
\newcommand{\eL}{{\mathrel{\mathscr L}}} 

\usepackage{xcolor}

\def\malce{\protect\mathbin{\hbox{\protect$\bigcirc$\rlap{\kern-7.75pt\raise0,50pt\hbox{\protect$\mathtt{m}$}}}}}

{\theoremstyle{definition}
}
{\theoremstyle{remark}
}

{\theoremstyle{remark}
}

\begin{document}

\sloppy

\title{The $q$-theory of Finite Semigroups: History and Mathematics}

\author{Stuart Margolis}
\address[S.~Margolis]{%
    Department of Mathematics\\
    Bar Ilan University\\
    52900 Ramat Gan\\
    Israel}
\email{margolis@math.biu.ac.il}

\maketitle

The book under review, {\em The q-theory of Finite Semigroups} by John Rhodes and Benjamin Steinberg is a remarkable achievement. It combines fresh new ideas that give a new light on 40 years of the study of pseudovarieties of finite semigroups and monoids together with the first comprehensive review of the decomposition and complexity theory of finite semigroups since 1976. The authors lay out a new approach that unifies and expands many threads that have run through the study of pseudovarieties and their applications. In addition, a number of fundamental results about finite semigroups that have never appeared in book form before are included and written from a current point of view. The letter ``$q$" stands for {\em quantum} and indeed, the main innovation in the book is to study continuous operators (in the sense of lattice theory) on the lattice of pseudovarieties and relational morphisms, giving an analogue to quantum theory. 

This paper is an expanded version of the short review, \cite{Shortqtheory}. A mathematics book does not appear independently, but within a historical and mathematical context. Before giving a detailed review of the contents of the book, I will give a short history of books in finite semigroup theory and on the theory itself in order to place this book in its proper context. There does not seem to be a traditional way to publish a book review this long. Thus, I am ``self publishing" this review. Email feedback to the author from readers would be greatly appreciated and I will update the paper when necessary.

\section*{Books on semigroup theory}

The first book in finite group theory is certainly {\em Trait\'{e} des substitutions et des \'{e}quations alg\'{e}briques (Treatise on Substitutions and Algebraic Equations)} by Camille Jordan \cite{Jordan} first published in 1870. Burnside's classic, {\em The Theory of Groups of a Finite Order} \cite{Burnside}, was first published in 1897 with its even more important second edition (because of its inclusion of the representation theory of finite groups for the first time) in 1907. Hundreds of books have been written with the phrase ``finite groups" in the title. These include books on the general theory, on representation theory by permutations and matrices, character theory, connections with finite geometry and other combinatorial structures as well as many other topics.

In contrast, there is no book with the title akin to that of Burnside's, simply ``{\em The Theory of Finite Semigroups}". As far as I know the only books in which the phrase ``finite semigroups" appears in the title are Almeida's {\em Finite semigroups and universal algebra}, \cite{Almeida:book},  Ganyushkin and  Mazorchuk's, {\em Classical finite transformation semigroups, an introduction} \cite{GM} and the current book. There are a number of books in which finite semigroups play a central role, all related to the interaction with finite semigroup theory and theoretical computer science, for example, \cite{Arbib, Eilenberg, Lallement, Pinbook, Straubingbook}.

While the monograph by Suschkewitsch \cite{Susch} is arguably the first book in which semigroups appeared as independent objects to be studied for their own sake, it was unavailable in the West. The first available books and monographs on algebraic semigroup theory were comprehensive. The books by Clifford and Preston and Lyapin educated Western and Soviet Block semigroup theorists for a generation \cite{CP, CP2, Lyapin1, Lyapin2}. In 1970, John Rhodes wrote a comprehensive review \cite{RhodesReferee} on four books that existed at that time. Besides the books by Clifford and Preston and Lyapin, reviews of Hofmann and Mostert's book on compact semigroups and R{\'e}dei's book on finitely generated commutative semigroups \cite{HofMost, Redei} were discussed. The reader would profit by reading Rhodes's review to get an impression on what was available in the late 1960s and how the previous generation of semigroup theorists was educated.

The first book dedicated to finite semigroups and their relationship to finite automata and formal languages was ostensibly a conference proceedings edited by Arbib \cite{Arbib}. However, chapters, 1, 5-9, written by John Rhodes and Bret Tilson, were a comprehensive introduction to finite semigroup theory. Through the mid-1970s, for students of John Rhodes, like the author of this review, this book was the main source of information on the subject. Chapter 4 of the book under review can be regarded as a major updating of much of the material in these chapters as well as Tilson's treatment of complexity in 1976 \cite{TilsonXI, TilsonXII}.

In the 1970s a few more books  were written on semigroup theory. Petrich's two books were dedicated to topics in the general structure theory of semigroups \cite{Petrich1, Petrich2}. Howie's book \cite{Howie1} was an updated version of Clifford and Preston, albeit less encyclopedic. None of these books dealt at all specifically with finite semigroups.

The major contribution in this era were the books by Eilenberg, especially Volume B with its two chapters by Tilson \cite{EilenbergA, Eilenberg, TilsonXI, TilsonXII}. These were originally planned as a four volume series on automata and formal language theory, but only the first two volumes appeared. The first volume is a reworking of the theory of finite automata, heavily influenced by the work of Sch\"{u}tzenberger and the French School. Volume B, heavily influenced by Sch\"{u}tzenberger and Tilson was dedicated to the theory of regular languages and Krohn-Rhodes theory.

The presentation in Volume B had many important innovations. One of the most important was the introduction of the notion of pseudovarieties of finite semigroups and finite monoids and their partner, the varieties of formal languages. The Eilenberg-Sch\"{u}tzenberger correspondence sets up a one-to-one correspondence between the algebraic and language theoretic varieties. The book gives the most important examples of that correspondence known at that time (besides Kleene's Theorem): Sch\"{u}tzenberger's Theorem relating aperiodic semigroups and star-free languages, and Simon's Theorems, one relating piecewise testable languages and $\mathcal{J}$-trivial monoids and one relating locally testable semigroups and local semilattices. It is not an overstatement to say that since 1976 with the appearance of Eilenberg's Volume B, the vast majority of finite semigroup theory has been involved with the study of pseudovarieties of finite semigroups and monoids and their relationship to automata theory.

Tilson's chapters \cite{TilsonXI, Tilson XII} were devoted to the Krohn-Rhodes complexity of finite semigroups. They gave a new proof of the fundamental lemma of complexity that introduced the Rhodes expansion, relational morphisms and the derived semigroup, which is a precursor of the derived and kernel categories. All of these have played an important role in finite semigroup theory ever since. Tilson also gave a notion of how a class of relational morphisms (weakly closed classes in his terminology) can be used to define operations on pseudovarieties.

The chief innovation of the book under review is to give a precise foundation for such classes of relational morphisms as well as determining the corresponding operators on the lattice of pseudovarieties. In addition, the book gives a complete update of classical Krohn-Rhodes complexity, the first since Tilson's chapters and also a generalization to so called two-sided complexity, more general complexity hierarchies and decomposition and complexity of finite idempotent additive semirings.

Another important book from this era is Lallement's {\em Semigroups and Combinatorial Applications} \cite{Lallement}. This book contains a standard introduction to semigroups, Green-Rees theory and Krohn-Rhodes theory. It covers material from Eilenberg Volume B and from the proposed, but never published Volume C, which gives an algebraic and combinatorial approach to context free languages. It also has semigroup theoretic applications to algebraic combinatorics.

From the 1980s onwards, many more books and conference proceedings have been published on semigroup theory. Of course, as the theory has matured, books tended to be on more specific topics rather than attempting a complete description of all of algebraic semigroup theory. It is not possible to discuss them all. Among the books that have influenced my own work are the following. Pin's book \cite{Pinbook} is an excellent introduction to the theory of varieties at an introductory graduate level. It is the first book that gives detailed descriptions of operations on varieties of one type and the corresponding operations on the second type of variety. Almeida \cite{Almeida:book} gives the first book treatment of the many deep and fundamental connections between pseudovarieties and profinite semigroups, as well as giving many more examples of operations on varieties. Putcha \cite{Putcha} and Renner \cite{Renner} develop the theory of linear algebraic monoids as well as their finite analogues. Okninski \cite{okninski} is devoted to matrix semigroups. Higgins \cite{higginsbook} includes lots of material that is not available in book form elsewhere. Grillet's book \cite{grillet} contains a chapter devoted to finite semigroups, including the only book version of the Rhodes-Allen synthesis theorem  \cite{RhodesAllen} as well as its version by Birget \cite{BirgetSyn1, BirgetSyn2}. Petrich \cite{PetrichInv} wrote an encyclopedic book on inverse semigroup theory. Lawson's book on inverse semigroups \cite{Lawson} emphasizes the connection to the theory of partial symmetries with numerous applications to many areas of mathematics. Petrich and Reilly's work on completely regular semigroups \cite{CRbook} details semigroups that are unions of their subgroups. Ganyushkin and  Mazorchuk's book \cite{GM} gives a modern treatment of naturally arising transformation semigroups and their representations of various kinds.

\section*{A short history of finite semigroup theory}

In this section, I give a short history of finite semigroup theory geared to explaining where the problems considered in the book under review came from. I assume that the reader has a basic knowledge of semigroup theory. Up until 1962, finite semigroup theory consisted of the Green-Rees structure theory and results on the representation theory of finite semigroups, due to Clifford, Munn and Ponizovsky. Although Clifford and Preston Volume 1 \cite{CP} was not specifically about finite semigroups, it contained essentially everything known about finite semigroups (modulo groups of course) at the time of its publication in 1961.

The Suschkewitsch Theorem \cite{Suschkewitsch} was perhaps the first theorem about the structure of finite semigroups. It determined the structure of the minimal ideal of a finite semigroup and in modern terminology, the structure of finite completely simple semigroups. The Rees Theorem \cite{Rees} and Green's relations \cite{Green}, although concerned with more general semigroups, are crucial for the understanding of finite semigroups. They determine the local structure of finite semigroups. Each finite semigroup $S$ has a principal ideal series and each principal factor corresponds to a unique $\J$-class of $S$ and is either a null semigroup or a completely 0-simple semigroup. The structure of the latter is determined by Rees Theorem.

Green-Rees theory tells one nothing about the global structure of a finite semigroup- that is, it doesn't tell one how elements from different $\J$-classes interact. The first advance from the local point of view were the Sch\"{u}tzenberger representations \cite{Schutzrep, Schutzmonomial}. These representations  show how elements $\J$-greater or equal to a fixed $\J$-class $J$ act on $J$. The importance of the Sch\"{u}tzenberger representations was not recognized in 1961. Indeed, Clifford and Preston wrote \cite{CP} on page 87: ``No real use has as yet been made of the Sch\"{u}tzenberger representation in semigroup theory." It was only recognized later that the Sch\"{u}tzenberger representations play an intimate role in describing the simple modules of a finite semigroup over a field, \cite{LallePet, myirreps}. Ironically, Chapter 5 of \cite{CP} is devoted to matrix representations of finite semigroups. Also, the associated Sch\"{u}tzenberger graphs and automata play the same role in semigroup theory and inverse semigroup theory as the Cayley graph does in group theory.

More closely related to the contents of the current book, the Sch\"{u}tzenberger representations are by monomial matrices with entries in the Sch\"{u}tzenberger group of a $\D$-class $D$ and thus in an appropriate wreath product of this group with the action of the semigroup on the $\eL$ or $\R$ classes in $D$. This gives a deep connection to Krohn-Rhodes theory as we will soon note. The use of the Sch\"{u}tzenberger representations is so important to Krohn-Rhodes theory, that it is called semilocal theory. Chapter 8 of \cite{Arbib} describes this and its uses in detail and the current book gives an updated description of this as part of Chapter 4.

1962 was a milestone in the history of finite semigroup theory with the appearance of the Krohn-Rhodes Theorem,
originally called the Prime Decomposition Theorem. It has had lasting influence on semigroup theory and automata and formal language theory
for over half a century and still provides the framework for current research.
While the connection between semigroups and sequential machine decompositions
had been looked at by a number of researchers in the late 1950s and early
1960s, the theorem established a definitive connection between machine
decompositions and wreath product decompositions of their syntactic semigroups.
Thus any result about decompositions of machines can be translated into
a theorem about wreath product decompositions of finite semigroups. Conversely,
any wreath product decomposition of finite semigroups has a corresponding interpretation
in terms of factorizations of finite state machines.

Furthermore, the Krohn-Rhodes Theorem shows that any finite semigroup $S$,
divides (that is, is an image of a subsemigroup of) a wreath product of simple
groups that themselves divide $S$ and the three element monoid $U_2$ consisting of 2
right zeroes and an identity. Lastly, the only finite semigroups that are prime with
respect to the division ordering and wreath product (that is, semigroups which
whenever they divide the wreath product of two other semigroups must necessarily
divide one of the factors) are precisely the finite simple groups and the divisors
(which turn out to be the subsemigroups) of $U_2$. This is the reason for the name
Prime Decomposition Theorem. The theorem had tremendous impact in mathematics in general and semigroup theory
in particular and its connections with
computer science. It appears, for example, in the Encyclopedia Britannica under
Automata Theory.

The notion of division was a key contribution of Krohn and Rhodes; previous decomposition theorems, based on taking subsemigroups, could not decompose the semigroup of all maps on a finite set. The Krohn-Rhodes theorem is a semigroup analogue of the Krasner-Kaluzhnin Theorem that states that every extension of a group $N$ by a group $G$ embeds in the wreath product $N \wr G$. As a consequence every finite group embeds into the wreath product of its Jordan-Holder factors. However, the semigroup theorem crucially depends on division as opposed to embedding.

The theorem leads to a definition that has occupied John Rhodes
in one form or another for the rest of his career. This is the {\em (group) complexity} of
a finite semigroup. It follows from the Krohn-Rhodes Theorem that every finite
semigroup divides a wreath product of groups and aperiodic semigroups, the latter
being semigroups with trivial maximal subgroups. The smallest number of groups
in any such decomposition is the (group) complexity of the semigroup. It was known
early on that there are semigroups of arbitrary complexity. For example, the full transformation semigroup
on $n$ elements has complexity $n-1$. However, since one needs to
search an infinite number of possible decompositions, there is no obvious way to
compute the complexity of a finite semigroup, given its multiplication table. This
problem has occupied Rhodes's work for 50 years, and many of his contributions have been directed by
the question of whether there is an algorithm to compute the complexity of a finite semigroup.

The search for the solution to the decidability problem, like many great mathematical
problems, has led to the development of many tools and ideas that are both
of direct use in finite semigroup theory and of independent interest. Thus the use
of techniques from category theory, topology, profinite algebra, logic, automata theory, ordered sets and
other fields have become important parts of finite semigroup theory. Understanding how and why these other
fields became necessary tools to study finite semigroups is important for understanding the problems studied in the current book.

A key innovation is the notion of a relational morphism between semigroups. A relational morphism $\phi\colon S \to T$ between semigroups $S$ and $T$ is a subsemigroup of $ S \times T$ that projects onto $S$. If we identify a usual homomorphism of semigroups with its graph, then we see that relational morphisms are generalizations of homomorphisms. Relational morphisms were first defined for groups by Wedderburn \cite{WedRelation}. They also appear in homological algebra (slightly generalized so that the projection to the first coordinate need not be onto) under the name of additive relation of modules \cite{MacLane}, Chapter 2.6, \cite{MacLane2} or correspondence \cite{Puppe}. The connecting homomorphism in the long exact sequence of homology can be defined as a composition of additive relations and many other basic features of homological algebra have an interpretation in terms of additive relations.

However, relational morphisms between groups or modules don't give too much that is new and are mainly a convenience. Thus, if $\phi\colon G \to H$ is a surjective relational morphism between groups, then $G/\phi^{-1}(1)$ is isomorphic to $H/\phi(1)$ \cite{WedRelation} with a similar relationship for additive relations between modules (and similar constructions in arbitrary Abelian categories). Thus in these categories, relational morphisms are essentially isomorphisms between divisors. In this form, the theorem goes back to Goursat's Theorem for groups from 1895 \cite{Lambek}. In fact, one can prove that a similar situation holds in any variety of universal algebras (where relational morphisms are subalgebras of the direct product) in which all congruences on all of its members commute \cite{Lambek}.

Relational morphisms between semigroups are much more complex to understand. They have played a central role in semigroup theory since their appearance in the mid-1970s and indeed the present book is intensely involved with relational morphisms. Where do they come from in finite semigroup theory? First of all, it is easy to see that a relational morphism $\phi\colon S \to T$ between semigroups is injective in the sense that $\phi(s_{1}) \cap \phi(s_{2}) \neq \emptyset$ implies that $s_{1} = s_{2}$ for all $s_{1}, s_{2} \in S$ if and only if $\phi$ is the inverse of a partial function and thus if and only if $\phi^{-1}:\phi(S) \to S$ is a surjective homomorphism. Thus $S$ divides $T$ if and only if there is an injective relational morphism $\phi\colon S \to T$. Therefore division of semigroups is naturally modeled by the class of injective relational morphisms.

In semigroup decomposition theory, relational morphisms arise naturally as follows. Suppose that $S \prec T\wr U$ is a division of $S$ into the wreath product of $T$ and $U$. Then there is an injective relational morphism $d: S \to T\wr U$. Let $\pi: T\wr U \to U$ be the projection. We thus obtain the relational morphism $\phi=\pi\circ d:S \to U$. The Krohn-Rhodes Theorem and complexity theory lead one naturally to study such relational morphisms when $T$ and $U$ are groups or aperiodic semigroups. The crucial question is the converse. Suppose that $\phi\colon S \to U$ is an arbitrary relational morphism. How can one find a ``minimal" $T$ to complete this to a wreath product division $S \prec T\wr U$?

If one looks at the case where $S$ is fixed and $U$ runs through all finite groups, then one is naturally led to the type II subsemigroup and eventually to deep connections with profinite group theory and the theorems of Ash and of Ribes-Zaleskii \cite{Ash, RZ}. If one looks at the case where $T$ is an aperiodic semigroup, we arrive to the Fundamental Lemma of complexity, the Rhodes expansion and other expansions, Sch\"{u}tzenberger's Theorem on star-free languages and generalizations \cite{Eilenberg, Tilson XII}. If $S$ is fixed and $U$ runs through aperiodic semigroups, we are led to the notion of pointlike aperiodic sets and Henckell's Theorem \cite{Henckell}. We eventually arrive to the derived category \cite{Tilson} and the kernel category \cite{Kernel} of a morphism and the introduction of methods from algebraic category theory into semigroup theory. Space doesn't allow me to expand on all of these topics, but all of this is described in detail in the book under review and is indeed one of the book's important contributions.

The Eilenberg-Sch\"{u}tzenberger Correspondence Theorem \cite{Eilenberg} showed that pseudovarieties (of semigroups, monoids and languages) played a central role in the theory. While pseudovarieties of semigroups and monoids are the finite version of the corresponding varieties in the sense of universal algebra, there was not a ready-made equational theory that parallels Birkhoff's Theorem \cite{universalalgebra}. The original version of Eilenberg and Sch\"{u}tzenberger was in terms of ultimately satisfying a sequence of identities $\{u_{i} = v_{i}\mid i \in I\}$, where $u_{i}$ and $v_{i}$ are elements of a countably generated free semigroup (or monoid). A semigroup ultimately satisfies such a sequence if it satisfies, in the usual sense, all but a finite number of them. This description is not very useful on both the theoretical and computational level. There is no obvious definition of finitely based pseudovariety in this framework, for example. This was rectified by Reiterman's Theorem, that showed that one has an equational theory if one replaces free discrete semigroups by free profinite semigroups. That is, pseudovarieties are exactly the classes of finite semigroups that can be defined as satisfying a set $\{u_{i} = v_{i}\mid i \in I\}$, where now $u_{i}$ and $v_{i}$ are elements of a free profinite semigroup \cite{Reiterman}. Of course, this is true in the more general context of pseudovarieties of finite universal algebras.

This may seem to be a pyrrhic victory in that describing elements of free profinite semigroups are certainly much more complex than those of free discrete semigroups.
However, Jorge Almeida and later his colleagues in a long sequence of papers showed that syntactic profinite methods were incredibly powerful and a necessary tool for anyone interested in pseudovarieties of finite semigroups. Much of this up to 1992 is summarized in Almeida's book \cite{Almeida:book}. Along with the connections to the type II connection mentioned above, this work made it clear that profinite methods were also a central part of the study of finite semigroups.

Finally, Jean-Eric Pin and coworkers, especially Pascal Weil, introduced methods of ordered semigroups and monoids into the mix \cite{PinOrd}. This is completely natural from the point of view of applications to language theory, as every syntactic monoid is ordered (in two dual ways) by splitting the definition of the syntactic congruence into two implications. This has allowed for generalization of the Eilenberg-Sch\"{u}tzenberger Theorem to classes of languages not necessarily closed under complementation.

Thus, in the generation Sch\"{u}tzenberger since the appearance of Eilenberg's Volume B \cite{Eilenberg}, methods from category theory, profinite algebra and ordered structures have become a class of absolutely essential tools for workers in finite semigroup theory and its applications. I will now try to set the stage for the problems that led to the central issues in the current book.

The Eilenberg-Sch\"{u}tzenberger Correspondence Theorem \cite{Eilenberg} arose from the examples mentioned previously: Kleene's Theorem that says that a language is regular if and only if its syntactic semigroup is finite, Sch\"{u}tzenberger that says that a language is star-free if and only if its syntactic semigroup is a finite and aperiodic semigroup and the two aforementioned theorems on piecewise testable languages and $\J$-trivial semigroups and locally testable languages and local semilattices. Thus one was led to ask what properties of a collection of rational languages that allowed it to be classified by its syntactic invariants. One definitive answer is the Sch\"{u}tzenberger Correspondence Theorem \cite{Eilenberg}. After the appearance of Eilenberg Volume B, \cite{Eilenberg} there were many papers written that gave examples of the correspondence between pseudovarieties of semigroups and monoids and varieties of rational languages.

Very soon a new direction was initiated. This was a study of how operations on one sort of variety was manifested on the variety of the second kind. Thus Straubing \cite{Straubing} proved that the concatenation of languages corresponded to the operation that sends a variety of semigroups {\bf V} to the variety ${\bf Ap \malce V}$, the Malcev product of the pseudovariety {\bf Ap} aperiodic semigroups and {\bf V}. From the point of view of the current book, this variety is defined as the collection of all finite semigroups $S$ for which there is a relational morphism $\phi\colon S \to T$ such that $T \in {\bf V}$ and $\{\phi^{-1}(e)\mid e=e^{2} \in T\} \subset {\bf Ap}$. More generally, one defines  the  Malcev product of two varieties ${\bf W \malce V}$ similarly.

Semidirect products of semigroups also lead to an operation on pseudovarieties of semigroups. Namely, if {\bf V} and {\bf W} are pseudovarieties, then the product {\bf V*W} is the collection of all finite semigroups that divide a semidirect product $S\rtimes T$, where $S \in {\bf V}$ and $T \in {\bf W}$. It can be proved that this operation is associative on the lattice of pseudovarieties of finite semigroups. If $V$ and $W$ are (pseudo)varieties of groups, then it is known (with use of the Krasner-Kaloujine Theorem) that ${\bf V \malce W}$ = {\bf V*W} and thus the Malcev product is an associative product on the lattice of (pseudo)varieties of (finite) groups.

It is to be expected the relationship between semidirect product and Malcev product is much more complex for pseudovarieties of semigroups. Trying to understand the
connections was fundamental to the development of the derived category and the kernel category. For example, the fundamental lemma of complexity \cite{TilsonXII} states that if {\bf V$_{n}$} is the pseudovariety of finite semigroups of complexity at most $n$, then ${\bf Ap \malce V_{n}}$ = ${\bf V_{n}}$. In fact, the derived category (in its
primitive form as the derived semigroup) and the Rhodes expansion \cite{TilsonXII} were invented in order to prove the equation
${\bf Ap \malce V_{n}} = {\bf Ap*(V_{n}*_{r}Ap)}$, where  $*_r$ is the reverse semidirect product operator on pseudovarieties of semigroups.

The derived category \cite{Tilson} gives the correct generalization of the relationship between Malcev product and semidirect product of group pseudovarieties to pseudovarieties of semigroups. The derived category theorem \cite{Tilson} shows that if {\bf V} and {\bf W} are pseudovarieties of semigroups, then
$V*W = \{S\mid \exists \phi\colon S \to T, \text{a relational morphism, such that} T \in {\bf W} \text{and the derived category} D_{\phi} \in {\bf gV}\}$. Here {\bf gV} is the pseudovariety of finite categories that divide members of $V$. For details on the derived category and the meaning of division of categories see \cite{Tilson} or the current book \cite{qtheor}. A similar relationship holds between the two-sided semidirect product operator on the lattice of pseudovarieties that sends the pair of pseudovarieties {\bf V,W} to {\bf V**W}, the pseudovariety of semigroups generated by two-sided semidirect products of members of {\bf V} and {\bf W} by replacing the derived category with the kernel category above \cite{Kernel, qtheor}.

The point of the previous paragraph is not to detail the properties of derived and the kernel categories, but to point out that the natural operations on the lattice of pseudovarieties associated to semidirect products are defined by extending one pseudovariety by a suitable collection of relational morphisms. This relationship is clear from the definition of the Malcev product. It also holds true for the join via Steinberg's notion of slice \cite{Slice}. Thus we finally arrive at the central problem studied in the current book \cite{qtheor}: to clarify the relationship between classes of relational morphisms and operators on the lattice of pseudovarieties of finite semigroups, with a theory that includes the examples cited here as special cases, as well as as many of the other situations that have arisen in the literature.

\section*{Review of The $q$-theory of finite semigroups}

Although the current book has four sections, it reads as an interweaved collection concerned with two main topics. Part I- ``The $q$-operator and Pseudovarieties of Relational Morphisms" and Part III- ``The Algebraic Lattice of Semigroup Pseudovarieties", deal with coming to grips with the problems detailed at the end of the previous section. The goal is to axiomatize classes of relational morphisms that define, via extension, nice operators on the lattice of pseudovarieties of semigroups. Part II- ``Complexity in Finite Semigroup Theory" and Part IV- ``Quantales, Idempotent Semirings, Matrix Algebras and the Triangular Product" give a more or less complete picture of the decomposition and complexity theory of finite semigroups up to the time of publication of the book and a generalization to finite semirings. Here the language and basic results of the newly developed $q$-theory are used to clarify some of the results and simplify some proofs. It is convenient to review each of these substrata separately.

\subsection*{$q$-theory and continuous operators on the algebraic lattice of pseudovarieties}

The main point of Part I and Part III is to study (Scott) continuous \cite{CL}  operators (in the sense of lattice theory) on the lattice {\bf PV} of pseudovarieties of finite semigroups. The operations mentioned above: semidirect product, two-sided product, Malcev product and join are certainly among the most intensively studied operators on {\bf PV} and all are continuous. Furthermore, in analogy with the Derived Category Theorem and The Kernel Category theorem, the definition of Malcev products and the Slice Theorem, well behaved collections of relational morphisms exist which define these operators by extension. The main question addressed in these parts of the book is to axiomatize the classes of relational morphisms that one needs to define all continuous operators. The emphasis on having continuous operators on {\bf PV} as the central objects of study gives an analogy to quantum mechanics. Indeed the mysterious ``$q$" in the book's title is a reference to this connection. Composition of operators is thus associative and one thinks of evaluating an operator at a pseudovariety as sampling. So one replaces points (= pseudovarieties) by operators and gets back points by sampling.

The most general definition relating relational morphisms to operators on the lattice of pseudovarieties along the lines suggested by the derived and kernel category
theorems, the Malcev product and the Slice Theorem is the following. Let $R$ be a collection of relational morphisms. Define $Rq$ to be the operator that takes a
pseudovariety {\bf V} to the collection $\{S\mid \exists \Phi:S \to T, \Phi \in R, T \in {\bf V}\}$. This is clearly taking the motivating examples and turning the
results into a definition.

Of course, in order to have a reasonable theory, one must put restrictions on the classes of relational morphisms and on operators. One would need some axioms to ensure that $Rq$ sends a pseudovariety to a pseudovariety and that $Rq$ had nice properties as an operator on {\bf PV}. The authors choose the Scott continuous functions on {\bf PV} \cite{CL}. Scott continuous functions between arbitrary posets have a rich and well developed  theory that has many applications in mathematics and computer science. Furthermore, naturally defined operators on {\bf PV} are continuous. 

The question then is how to define classes of relational morphisms ${R}$ such that their corresponding operator ${Rq}:{\bf PV} \to {\bf PV}$ is continuous. Then one would want to know exactly what operators are defined this way and study how this operation effects pseudoidentities defining {\bf V} and $Rq(\bf V)$. These questions among others comprise Part I of the book.

Chapter 1, ``Foundations for Finite Semigroup Theory" is concerned with background material in categories, wreath and one and two-sided semidirect products and basic ideas from order theory. Green-Rees theory is relegated to Appendix A, so it would help a reader to already have a foundation in classical semigroup theory from other sources before reading Chapter 1. Chapter 1 is concerned with setting up the categorical framework for studying the category whose objects are finite semigroups and whose morphisms are relational morphisms. When dealing with the corresponding category of monoids and monoid relational morphisms, one must be careful to use the adjoint of the forgetful functor from monoids to semigroups, which amounts to adding an external identity to a semigroup even if it already has an identity. Using the more familiar ``add an identity to a semigroup only when it doesn't have one already" has led to a number of subtle errors in published proofs and definitions. Thus category theory is considered to be the correct formulation for moving between categories and is done so via adjoint functors when they exist.

The treatment of relational morphisms is crucial to the book. A notion of division between relational morphisms and relational morphisms between relational morphisms is defined. Thus the category {\bf (FSgps, RM)} of finite semigroups and relational morphisms, becomes a 2-category in that the collection of morphisms between any two objects is itself a category. As one might expect, the category {\bf (FSgps, RM)} is more complex to understand than the classical category {\bf (FSgps, FM)}, where {\bf FM} is the collection of functional homomorphisms.

This sets the stage for Chapter 2, ``The $q$-operator". This chapter defines the classes of relational morphisms that allow one to define continuous operators on {\bf PV} as alluded to above. These are defined as satisfying naturally defined axioms that mimic the operations defining pseudovarieties of semigroups. A {\em pseudovariety of relational morphisms} is a class of relational morphisms closed under product of relational morphisms, division of relational morphisms and range extension. The first two operations are described in detail in Chapter 1 and the last just says that if $\phi\colon S \to T$ is in a class of relational morphisms and $j:T \to T'$ is a functional injection, then $j\phi\colon S \to T'$ is also in the class. A {\em continuously closed class of relational morphisms} is defined similarly, except that a stronger form of division is required.

Now the $q$-operator is defined for these classes as above. That is, if $R$ is a pseudovariety of relational morphisms or a continuously closed class, then
$Rq:{\bf PV} \to {\bf PV}$ is the operator $\{S\mid \exists \Phi:S \to T, \Phi \in R, T \in {\bf V}\}$. One proves that for pseudovarieties and continuously closed classes of relational morphisms, the image of a pseudovariety is a pseudovariety. Moreover, $Rq$ is a continuous operator on ${\bf PV}$. Most importantly, every continuous operator on ${\bf PV}$ is of the form $Rq$ for some continuously closed class $R$. The image of the collection of pseudovarieties of relational morphisms is the collection of operators satisfying the so called generalized Malcev condition. These form a proper collection of operators, but include the ones usually studied in the literature such as semidirect products, two-sided products and generalized Malcev products. An important example of an operator that does not satisfy the generalized Malcev condition is the well studied operator that sends a pseudovariety {\bf V} to the pseudovariety {\bf PowV} generated by the power semigroups of members of {\bf V}.

The collection of pseudovarieties of relational morphisms {\bf PVRM}, continuous closed classes of relational morphisms {\bf CC}, continuous operators {\bf Cnt(PV)} and {\bf GMC(PV)} of continuous operators that satisfy the generalized Malcev condition all form algebraic lattices. The $q$-function is order preserving. An important property is that $q$ is the left part of a Galois correspondence between these lattices. This allows one to prove that the image of the $q$ function is all of {\bf Cnt(PV)} if the domain is {\bf CC} and all of {\bf GMC(PV)} if the domain is {\bf PVRM}. It is shown that $q$ is not injective, that is different collections of relational morphisms can define the same continuous operator. Being part of a Galois correspondence proves that there are minimal and maximal such collections mapping onto the same operator that are of natural interest. Computing them is highly non-trivial, even for the smallest class of relational morphisms, which is the collection of divisions, whose identification as the minimal representative of the identity operator requires a non-trivial theorem of Bergman \cite{Bergman}. A word of warning- the determined join in the lattices of operators is not the infinite intersection. This has lead to errors in the literature. A description of the determined join is given.

Chapter 2, ``The $q$-operator" continues with a collection of important examples of the various lattices above. The derived category theorem and the kernel category theorem are proved, and are generalized to semigroupoids, that is, categories that may not have identities at each object. These are the motivating examples for $q$-theory. There is a natural associative multiplication on {\bf PVRM} and {\bf CC} induced by composition of relations and $q$ is then a monoid morphism onto {\bf Cnt(PV)} and {\bf GMC(PV)} respectively. The $q$ function can be modeled by composition and this is an important technical tool.

Chapter 3 is entitled ``The Equational Theory". As outlined in the history section of this review, Reiterman's Theorem gives the analogue of Birkhoff's Theorem for pseudovarieties. Pseudovarieties are defined by identities in free profinite semigroups. It is then natural to consider the following problem. Given a description of pseudoidentities defining pseudovarieties {\bf V} and {\bf W} and an ``operation" $\Box$ forming a new pseudovariety ${\bf V}\Box {\bf W}$, give a pseudoidentity basis for ${\bf V}\Box {\bf W}$. Of course, this can be generalized to $n$-ary operations on pseudovarieties as well. The motivating examples, as usual, are semidirect and Malcev products and indeed a lot of work over the last years have been devoted to finding pseudovariety bases for these operations \cite{Almeida&Steinberg:2000, Almeida&Steinberg:1998, AlmeidaWeil, PinWeil}. One hope for deciding complexity is to show that if one starts with aperiodic semigroups and groups, then decidable pseudovariety bases can be given for the complexity hierarchies by iterating bases relative to semidirect products and/or Malcev products. The purpose of this chapter is to show that these individual results are special cases of a general basis theorem for the composition of two pseudovarieties of relational morphisms. Thus the results developed in the previous chapter apply.

The first section of this chapter is an introduction to profinite topology and in particular to profinite semigroup theory. The second section gives a proof of Reiterman's Theorem. The rest of the chapter is devoted to generalizing Reiterman's Theorem to the case of pseudovarieties of relational morphisms. This involves the development of new material on profinite relational morphisms and on relational pseudoidentities. A successful Reiterman Theorem is proved for pseudovarieties of relational morphisms in section 3.5. A related result, but not as complete, is proved for continuous closed classes in the last section of the chapter.

The main purpose of this chapter is to use all of this material to prove a general basis theorem that has as special cases, bases theorems for semidirect products, Malcev products and joins. The statement and proof involve preliminary material on inevitable substitutions, an idea that goes back to Ash's original solution of the type II conjecture \cite{Ash}. Rhodes' notion of flows and their consequences is needed to give the basis theorem for semidirect products.

While lattice theory and in particular the notions of algebraic and continuous lattices were used in Chapter 2, Chapter 6, ``Algebraic Lattices, Continuous Lattices and
Closure Operators" is essentially a short course in these subjects written with the needs of the semigroup theorist, and in particular, readers of the present book in mind. The seventh chapter is entitled ``The Abstract Spectral Theory of {\bf PV}". By the abstract spectral theory of a complete lattice $L$, the authors mean the study of the space $\text{Spec}(L)$ of finite meet indecomposable elements of $L$ (except the top) endowed with the hull-kernel topology whose closed sets are of the form $V(l) = \{m \in \text{Spec}(L) \mid m \geq l\}$ with $l \in L$. The authors are particularly interested in the lattice of pseudovarieties of finite semigroups, the lattice of pseudovarieties of relational morphisms of finite semigroups and the lattice of continuous operators on the lattice of pseudovarieties of finite semigroups.

The authors study various forms of indecomposable elements in these lattices, (meet and join irreducibles, as well as their finite and strong counterparts). They complete the study of the important class of Kov\'{a}cs-Newman semigroups, a class of finite join irreducible semigroups. Irreducibility with respect to semidirect product is also discussed.

The lattices that are studied in this book have associative operations that are compatible with the lattice operations. This leads to the study of these objects within the context of quantales. This is the subject of Chapter 8, ``Quantales". Classically, a quantale is a complete lattice that is also a semigroup and such that all left and right translations preserve arbitrary joins. They form a model for pointless topology, in that the open sets of a topological space form a quantale with respect to intersection as multiplication. They arise naturally on the border of semigroup and semiring theory. The semiring of all subsets of any semigroup is a quantale, where union and intersection are the lattice operations. In particular, the semigroup of all languages over a finite alphabet is a quantale and this gives connections with language theory that are discussed in this chapter.

The authors need a weaker condition, in that they require that left and right translations be continuous, that is, preserve directed joins. With this definition, {\bf PV, PVRM, Cnt(PV)} and other lattices that are intensively studied in this book become quantales with respect to naturally defined multiplications. The authors review the basics of quantale theory. They then interpret some of the maps between these lattices and in particular the $q$ map within the context of quantale theory.

Stone duality in the semigroup context shows that the category of profinite semigroups (monoids) is dually equivalent to the category of (counital) Boolean bialgebras. The comultiplication is explicitly constructed. The collection of recognizable languages over a finite alphabet is given the structure of a bialgebra and the authors give an interpretation of the Eilenberg-Sch\"{u}tzenberger Correspondence Theorem within the context of bialgebras and suggest this as a topic of further research.

\subsection*{Complexity of Semigroups and Continuous Operators}

Chapters 4, 5 and 9 are devoted to decomposition and complexity theory in the sense of the Krohn-Rhodes Theorem. The classical Krohn-Rhodes complexity of a finite
semigroup $S$ is, as noted above, the least number of groups needed in any wreath product decomposition of $S$ in terms of groups and aperiodic semigroups. Since the
pseudovarieties {\bf $V_n$}, $n \geq 0$, of semigroups of complexity at most $n$ can be defined by iterating the continuous operators that sends a pseudovariety {\bf V}
to {\bf Ap}*{\bf V} and {\bf Gp}*{\bf V}, where {\bf Ap} ({\bf Gp}) is the pseudovariety of aperiodic finite semigroups (finite groups), this opens up the way to look at complexity hierarchies defined by iterating collections of arbitrary continuous operators. In short, Chapter 4 is devoted to the classical one-sided theory, Chapter 5 to the two-sided theory and more general hierarchies and Chapter 9 is a beginning of the decomposition and complexity theory of finite semirings.

Chapter 4, ``The Complexity Theory of Finite Semigroups" is  a mini-book, being over 170 pages long and can be used (modulo some use of material from the previous chapters) as an introduction to this field. It is indeed the first book version of this theory since Tilson's chapter \cite{TilsonXII} from 1976 and the most extensive work on Krohn-Rhodes complexity since Chapters 1, 5-9 of  \cite{Arbib}. Starting with a proof of the Krohn-Rhodes Theorem, this chapter includes some of the most important results of the last 40 years that have never appeared in book form. Complexity heirarchies are defined in general and in particular the Krohn-Rhodes complexity of finite semigroups is defined. The Krohn-Rhodes Theorem is used to prove Sch\"{u}tzenberger's Theorem on star-free languages and aperiodic semigroups. Stiffler's determination of the pseudovarieties generated by the semidirect product closure of the two element prime semigroups is also proved.

As mentioned in the historical section, the semilocal theory of a finite semigroup is centered around the Sch\"{u}tzenberger representations and in general, how the $\mathcal{J}$-classes interact. Section 4.6 is a very long one that contains a treasure trove of material that are among the most useful and important tools for studying finite semigroups. In section 4.7, the authors use this material to help elucidate the structure of subdirectly irreducible finite semigroups. These two sections are an updated version of the material in Chapters 7 and 8 of \cite{Arbib}.

Chapter 9 of \cite{Arbib} is devoted to describing a number of characterizations of the complexity of completely regular semigroups, all of which are computable. Completely regular semigroups are precisely the Malcev product of completely simple semigroups with semilattices by a famous result of Clifford \cite{CP}, and thus are exactly the regular semigroups that have a morphism onto a semilattice that separates $\mathcal{J}$-classes. Early on Rhodes had proved a more general result that showed that some of the results of Chapter 9 of \cite{Arbib} allow one to show that the complexity of a finite semigroup that has a functional morphism onto an aperiodic semigroup that separates regular $\mathcal{J}$-classes has decidable complexity. Semilocal theory shows that the collection of all such semigroups is the pseudovariety ${\bf LG \malce Ap}$, where {\bf LG} is the pseudovariety of local groups, that is, the collection of finite semigroups such that each submonoid is a group. Section 4.10 gives a proof of this theorem. The proof uses semilocal theory and a formal manipulation of certain Malcev products and is simpler than previously proofs.

Section 4.13 gives the first book version of Graham's Theorem \cite{Graham} on the structure of the idempotent generated subsemigroup of a completely 0-simple semigroup. The approach here follows the topological approach given by Houghton \cite{houghton}, who proved this theorem independently. In fact, this is a ``frequently rediscovered theorem" whose proof has appeared a number of times in the literature.

Lower bounds for complexity were first considered in the 1970s by Rhodes and Tilson \cite{lowerbounds1, lowerbounds2}. This lead to the notion of type I and type II
subsemigroups of a semigroup $S$. The type II subsemigroup of $S$, also known as its group kernel, consists of elements that are related to the identity in any relational morphism from $S$ to a finite group. Less well known are the type I subsemigroups of a finite semigroup $S$. These are the analogue of the type II subsemigroup with respect to relational morphisms to aperiodic semigroups. A function from the class of finite semigroups to the natural numbers that satisfies all the axioms of complexity theory, plus the axiom that states that the value of $S$ is the maximal of the values of the submonoids of $S$ is called a local complexity function. Rhodes and Tilson \cite{localcomplex}, proved that the largest local complexity function is the maximal number of type I subsemigroups in a chain of subsemigroups alternating absolute type I subsemigroups and non-aperiodic type II subsemigroups. A semigroup is absolute type I if it is a type I subsemigroup of itself. Section 4.12 contains a proof of this lower bound. It is proved as well that  it is decidable whether a semigroup is an absolute type I semigroup, a theorem first appearing in \cite{HMPR} and thus, along with Ash's Theorem \cite{Ash}, proves that the largest local complexity function is computable.

Clearly, the local complexity of a finite semigroup is a lower bound to its complexity. An example is given in Section 4.14 that shows that this bound is strict and Section 4.16 shows that the complexity of a finite semigroup can differ arbitrarily from its local complexity. Thus complexity is highly global and depends not just on the local submonoids $eSe$ of a finite semigroup $S$, but how they fit together within all of $S$.

These proofs use the Presentation Lemma. This is a tool invented by Rhodes in the 1970s and is detailed in Section 4.14, using a formulation of Steinberg. The Presentation Lemma has some of the deepest applications to finite semigroups, notably allowing for the construction of highly non-trivial examples of semigroups. Besides giving the non-locality of complexity, the Presentation Lemma is used in section 4.15 to prove Tilson's Theorem \cite{2J} proving that the complexity of a semigroup with at most 2 non-zero $\mathcal{J}$-classes is decidable.

Sections 4.17 and 4.18 are dedicated to the proof of the Ribes-Zaleskii theorem \cite{RZ} and its connection to Ash's Theorem \cite{Ash}, which confirmed that the type II conjecture is true. The type II conjecture stated that the type II subsemigroup of a finite semigroup $S$ is the smallest subsemigroup $T$ of $S$ containing the idempotents and such that if $sts=s, s,t \in S$, then $sTt \cup tTs \subseteq T$. Rhodes and Tilson proved that every regular element of the type II subsemigroup is in this subsemigroup \cite{lowerbounds2}. Pin and Reutenauer \cite{Pintop, PR} discovered a connection between the type II conjecture and a question on the profinite topology of finitely generated free groups. This asked whether the set product of a finite number of finitely generated subgroups of a free group is closed in its profinite topology. Pin and Reutenauer  proved that the type II conjecture follows from this assertion. The Ribes-Zaleskii Theorem proved this assertion on the free group. The proof here involves an excursion into combinatorial group theory, especially around the notion of Stallings foldings \cite{Stallings}. The proof follows that of Auinger and Steinberg \cite{KarlBenRZ}.

Let {\bf V} be a pseudovariety of semigroups and $S$ a semigroup. A subset $X$ of $S$ is {\it pointlike with respect to {\bf V}} if for every relational morphism $\phi: S \to T$ from $S$ to a member $T$ of {\bf V}, there is a $t \in T$ such that $X \subseteq t\phi^{-1}$. The connection of this notion with the type II conjecture is clear, and indeed, Ash's Theorem can be used to determine the pointlike sets with respect to the pseudovariety of finite groups. Much earlier, Henckell \cite{Henckell} determined the pointlike sets with respect to the pseudovariety of aperiodic semigroups. The purpose of the last section of Chapter 4 is to prove Henckell's Theorem. The proof is long and highly technical and follows the more recent \cite{ouraperioidcpointlikes} using profinite methods.

Chapter 5 begins with a general framework for complexity hierarchies obtained by iterating operators on the lattice {\bf PV}. This generalizes the classical Krohn-Rhodes complexity. Tilson considered a generalization of $p$-length from group theory to semigroup theory in his thesis \cite{folleyT} that fits into this framework.

The main purpose in this chapter is to develop a two-sided Krohn-Rhodes complexity theory. Classical complexity is inherently one-sided by virtue of its definition based on one-sided semidirect or wreath products. There are examples of semigroups $S_n, n > 0$, such that $S_n$ has complexity $n$, but its reverse has complexity 1. An example of such a sequence of semigroups is given in Chapter 4. Thus there is a need for a symmetric theory. One can use the two-sided semidirect product {\bf **} and iterate the operators {\bf Ap ** ( )} and {\bf Gp ** ( )}, but this runs into the problem that the two-sided semidirect product is not associative, so these operators are not idempotent. They do lead to the correct theory, but one must first iterate each of them to obtain their idempotent power (in the limit). This is a bit difficult to work with.

The goal of Sections 5.2 and 5.3 is to prove that these idempotent powers have simple descriptions in terms of (left) Malcev product operators associated with the pseudovarieties of aperiodic semigroups and local groups. The proof goes through part of the classification of maximal proper surmorphisms, that is, of minimal congruences of finite semigroups, first considered by Rhodes in \cite{MPS}. The papers \cite{RhodesWeil}and \cite{Kernel} influence the presentation in this book. This is material of independent interest. The proofs use much of the previous material, in particular, a heavy use of the semilocal theory and the kernel theorem. This justifies the explicit definition of the two-sided complexity hierarchy as the sequence of pseudovarieties $\{{\bf K_{n}}\mid n \geq 0\}$ as follows: ${\bf K_{0} = Ap, K_{n+1} = Ap \malce (LG \malce K_{n})}$. It should be noted that two-sided complexity has had many important applications in formal language theory. See, for example, \cite{Straubingbook, Unambig, Weil, Weilcounter}.

The rest of the chapter is concerned with showing that the two-sided complexity hierarchy is infinite. This again requires a non-trivial and computable lower bound theorem. First note that the two-sided complexity of a finite semigroup is at most its one-sided complexity. However, the two-sided complexity can be arbitrarily smaller than the corresponding one-sided complexity. For example, it is well known that the important pseudovariety {\bf DS} consisting of all finite semigroups whose regular $\mathcal{D}$-classes are subsemigroups satisfies the equation ${\bf DS = LG \malce Semilattices}$. Therefore, by definition every semigroup in {\bf DS} and in particular every finite completely regular semigroup has two-sided complexity at most 1. It is known that there are completely regular semigroups with arbitrary one-sided complexity.

A more difficult, previously unpublished result of Rhodes, is that every finite semigroup with at most 2 non-zero $\mathcal{J}$-classes has two-sided complexity at most 1. A lower bound result is proved and used to show that the full transformation monoid on $n$ elements has two-sided complexity $\lceil \frac{n-1}{2}\rceil$, proving that the two-sided hierarchy is indeed infinite. This is also an unpublished result of Rhodes. Despite its many applications in computer science, two-sided complexity theory is far less developed than the one-sided theory. Many open problems and challenges are given at the end of this chapter.

The last chapter of the book, ``The Triangular Product and Decomposition
Results for Semirings" is the beginning of a decomposition and complexity theory for finite idempotent semirings. A semiring is deemed idempotent if its additive monoid is idempotent and thus a lattice, if the semiring is finite. Finite idempotent semirings are exactly the finite quantales in the classical sense and this gives a connection to Chapter 8. The presentation in this chapter is from the point of view of finite idempotent semirings.

One needs an analogue of the wreath product in order to have a decomposition theory for semirings. Motivation comes from the classical theory of finite dimensional (associative) algebras over a field. If $A$ is such an algebra then any matrix representation over a field can be put into block triangular form, with irreducible blocks on the diagonal, by using the Jordan-H\"{o}lder Theorem for $A$ modules. This lead Plotkin and his coworkers to the definition and the study of the triangular product of $A$ modules \cite{Plotkin, Plotkinbook}. The definition of the triangular product works equally well over semirings and this is the product used for decomposition of idempotent semirings in this chapter.

Basic properties of semirings and triangular products are given. The triangular decomposition theorem of Section 9.3 shows that if $k$ is a finite idempotent semiring and $S$ a finite semigroup, then the semiring algebra $kS$ divides an $n$-fold triangular product of matrix semirings over the semiring of Sch\"{u}tzenberger groups of $S$. This is an exact analogue to the Munn-Ponizovsky Theorem for finite semigroup algebras over a field \cite{CP} and indeed the proof is modelled on the classical proof of this theorem via principal series of $S$ and Sch\"{u}tzenberger representations.

In Section 9.4 the results are specialized to the case that algebras are taken over the two element Boolean algebra $\mathbb{B}$, considered as a finite semiring. Notice that if $S$ is a finite semigroup, then its $\mathbb{B}$ algebra $\mathbb{B}(S)$ is isomorphic to the powersemiring $P(S)$. Powersemigroups have been extensively studied in finite semigroup theory, both for their algebraic properties and applications to formal language theory. See \cite{Almeida:book}. This shows that any finite $\mathbb{B}$ semiring divides an $n$-fold triangular product of matrix semirings over the powersemirings $P(G)$, where $G$ is a Sch\"{u}tzenberger group of the multiplicative semiring of $S$. The authors call this the Prime Decomposition Theorem for Finite Idempotent Semirings.

The authors study irreducibility of semirings with respect to triangular product. While they do not give a complete classification, it is proved that matrix semirings over the power semiring of a finite group are irreducible with respect to triangular product. In particular, this proves that like the case of finite groups and finite dimensional algebras over a field, but unlike the case of finite semigroups, every finite idempotent semiring embeds into an irreducible such semiring. Despite this, the authors give a definition of complexity for finite idempotent semirings. They use this to improve unpublished results of Fox and Rhodes \cite{FoxRhodes} on the Krohn-Rhodes complexity of the power semigroup of the powersemigroup $P(S)$ of a finite semigroup $S$. In particular, they give the exact answer in the case that $S$ is an inverse semigroup and give tight asymptotic bounds on the Krohn-Rhodes complexity of the powersemigroup of the full transformation semigroup on a finite set.

Semirings are clearly important in computer science. The natural numbers and the two element Boolean semiring come to mind immediately. In the last years, there has been an intensive study of the tropical semiring, consisting of the real numbers with $max$ as addition and the usual addition as multiplication. The work of Polak \cite{Polak1, Polak2} gives connections between finite idempotent semirings with formal language theory. It is thus natural to study representations of finite semigroups and finite idempotent semirings and this chapter gives an excellent direction for researchers to pursue.

Finally, the book contains a plethora of exercises, not all routine, open problems and historical and mathematical discussions that are useful for students and researchers alike. This work is a major contribution to the literature in finite semigroup theory and its applications. Both as a reference and as a research monograph pointing to new directions in research, it should be read in a deep way by anyone interested in this field.

\section*{Conclusion}

When I began my study of finite semigroup theory in 1975, I was able to learn all the background material I needed by reading the first volume of Clifford and Preston \cite{CP} and Arbib \cite{Arbib}. Of course, it was a big help to listen to the inspiring lectures of John Rhodes and have the wise counsel of the late Bret Tilson, where I learned about the then completely new ideas and uses of relational morphism and the derived semigroup (the precursor of the derived and kernel categories) in finite semigroup theory. The appearance of Eilenberg's Volume B with its two chapters by Tilson \cite{Eilenberg, TilsonXI, TilsonXII} in 1976 provided the necessary background and motivation to study pseudovarieties and pointed the way to research for my generation of finite semigroup theorists. The point here is that it was relatively easy for a beginning graduate student to learn the background material necessary to start doing research in finite semigroup theory.

This book is a testament to the amount of work that has been done in the past 40 years in this subject. In order to read this book in its entirety, a beginning graduate student (or any other reader) would have to have a background, not only in classical finite semigroup theory and its relationship to formal language theory and automata theory, but also a serious understanding of material in category theory, ordered sets and lattices, profinite topology and other areas of mathematics and computer science. It is certainly worth the effort to obtain the background to read this essential book for the study of finite semigroups.

This article is meant to put the work into its proper historical and mathematical context. The introduction of the $q$-operator is the authors' formalization of the connection between relational morphisms and operators on the lattice of pseudovarieties. This question goes back to Tilson's work \cite{TilsonXII}. The current book gives a deeply thought out theory that introduces many new ideas and results into finite semigroup theory. In addition, the book gives a comprehensive and up to date introduction to Krohn-Rhodes and complexity theory of finite semigroups, both to the classical theory and a foray into two-sided complexity and an introduction to decomposition and complexity theory of finite semirings. {\em The q-theory of Finite Semigroups} by John Rhodes and Benjamin Steinberg is a major achievement that will certainly have an influence on finite semigroup theory for the next generation of researchers.

{\bf Acknowledgement }I would like to thank Mark Lawson for his attention to this paper which led to many improvements.
\bibliography{stubib}
\bibliographystyle{amsplain}
\end{document}